\documentclass[10pt,conference]{ieeeconf}
\IEEEoverridecommandlockouts
\usepackage{amsmath,amssymb,latexsym,epsfig,psfrag,graphicx}
\newtheorem{thm}{Theorem}

\newenvironment{pf*}[1]{\smallbreak\noindent{\it #1}}{\phantom{.}\hfill$\Box$\smallbreak}
\newcounter{definition}
\newenvironment{dfn}{\addtocounter{definition}{1}\smallbreak\noindent
  {\em Definition \thedefinition.}}{\smallbreak}
\newcounter{remark}
\newenvironment{rmk}{\addtocounter{remark}{1}\smallbreak\noindent
  {\em Remark \theremark.}}{\smallbreak}
\newcounter{example}
\newenvironment{ex}{\addtocounter{example}{1}\smallbreak\noindent
  {\bf Example \theexample.}}{\smallbreak}
\newenvironment{ex*}[1]{\addtocounter{example}{1}\smallbreak\noindent
  {\bf Example \theexample{} --- {\bf #1}}}{\hfill$\Box$\smallbreak}
\newcommand{\realR}{\mathbb{R}}
\newcommand{\complexC}{\mathbb{C}}
\newcommand{\itemi}{($i$) }
\newcommand{\ii}{($ii$) }
\newcommand{\iii}{($iii$) }
\newcommand{\iv}{($iv$) }
\newcommand{\itemv}{($v$) }
\newcommand{\vi}{($vi$) }

\begin{document}
\title{Realizability and Internal Model Control on Networks}
  \author{Anders Rantzer
  \thanks{The author is affiliated with Automatic Control LTH, Lund
    University, Box 118, SE-221 00 Lund, Sweden.}}
\maketitle
\begin{abstract}
It is proved that network realizability of controllers can be enforced without conservatism using convex constraints on the closed loop transfer function.
Once a network realizable closed loop transfer matrix has been found, a corresponding controller can be implemented using a network structured version of Internal Model Control.
\end{abstract}

\section{Introduction}
The importance of closed loop convexity in the theory for control system design has long been recognized \cite{Boyd/B91}. A large number of specifications, both in time and frequency domain, can be stated as convex constraints on the closed loop system. The convexity opens up for efficient synthesis algorithms, as well as for computation of rigorous bounds on achievable performance. However, there are also many important specifications that cannot be expressed in a closed loop convex manner. A notable example is controller complexity, measured by the number of states needed in the realization. 

During the past decade, growing attention has been paid to large scale networks and distributed control. In this context, it is common to consider controller transfer matrices with a pre-specified sparsity pattern. Such sparsity constraints are generally not closed loop convex, but a number of important closed loop convex structural constraints have been derived and summarized under the framework known as quadratic invariance \cite{Rotkowitz06}. However, with exception for positive systems \cite{rantzer2018tutorial}, solutions based on sparsity restricted transfer matrices have a tendency to become computationally expensive and poorly scalable.
It was therefore an important discovery in the theory for large-scale control when \cite{wang2016system} recently proved that optimization with finite impulse response constraints can be used for scalable synthesis of distributed controllers. 

The objective of this short note is to isolate an idea used in ``system level synthesis'' \cite{wang2016system,anderson2017structured}, to show that network realizable controllers in the sense of \cite{vamsi2016optimal} can be synthesized using convex optimization. Moreover, we will demonstrate that the classical idea of internal model control \cite{garcia1982internal} is useful to convert the optimization outcome into a network compatible controller realization.


\section{Notation}
A \emph{transfer matrix} denotes a matrix of rational functions that can be written on the form $\mathbf{G}(z)=C(zI-A)^{-1}B+D$, where $A\in\realR^{n\times n}$, $B\in\realR^{p\times n}$, $C\in\realR^{n\times m}$ and $D\in\realR^{p\times m}$. It is said to be \emph{strictly proper} if $D=0$. 

\section{Network Realizability}

Following \cite{vamsi2016optimal}, we make the following definition.
\begin{dfn}
  Given a graph $\mathcal{G}=(\mathcal{V},\mathcal{E})$ with $N$ nodes, a transfer matrix $\mathbf{G}$ is said to be \emph{network realizable on} $\mathcal{G}$ if it has a stabilizable and detectable realization 
  \begin{align}
    \left[\begin{array}{c|c}
      A&B\\\hline C&D
    \end{array}\right]
    &=\left[\begin{array}{ccc|ccc}
      A_{11}&\ldots&A_{1N}&B_{1}&&0\\
      \vdots&&\vdots&&\ddots\\
      A_{N1}&\ldots&A_{NN}&0&&B_{N}\\\hline
      C_{11}&\ldots&C_{1n}&D_{1}&&0\\
      \vdots&&\vdots&&\ddots\\
      C_{N1}&\ldots&C_{NN}&0&&D_{N}
    \end{array}\right]
  \label{eqn:ABCD}
  \end{align}
  where $A_{ij}=0$ and $C_{ij}=0$ for $(i,j)\not\in\mathcal{E}$. Such a realization is said to be \emph{compatible with} $\mathcal{G}$. We need  $A_{ij}\in\realR^{n_i\times n_j}$, $B_i\in\realR^{n_i\times m_i}$, $C_{ij}\in\realR^{m_i\times n_j}$ and $D_i\in\realR^{p_i\times m_i}$, where $n_i$, $p_i$ and $m_i$ are the number of states, outputs and inputs in node $i$ respectively. 
\end{dfn}

Given a transfer matrix and a graph, no simple test for network realizability is known. However, given a realization, it is of course straightforward to verify the conditions of Definition~1.

\begin{thm}
Let the transfer matrices $\mathbf{G}_1$ and $\mathbf{G}_2$ be network realizable on $\mathcal{G}$. Then the following statements hold:
\begin{description}
  \item[\itemi] $\mathbf{G}_1+\mathbf{G}_2$ is network realizable on $\mathcal{G}$.
  \item[\ii] If $\mathbf{G}_1$ and $\mathbf{G}_2$ are stable, then $\mathbf{G}_1\mathbf{G}_2$ is network realizable on $\mathcal{G}$.
  \item[\iii] If $\mathbf{G}_1(\infty)$ is invertible, then $\mathbf{G}_1^{-1}$ is network realizable on $\mathcal{G}$.
\end{description}
\label{thm:add}
\end{thm}

\begin{rmk}
  Consider the graph with $\mathcal{V}=\{1,2,3,4\}$ and $\mathcal{E}=\{(1,1),(2,2),(3,3),(4,4),(1,3),(1,4),(2,3),(2,4)\}$. Notice that 
  both the two transfer matrices
  \begin{align*}
    \mathbf{G}_1(z)&=
    \begin{bmatrix}
      0&0&0&0\\
      0&0&0&0\\
      \frac{1}{z-2}&1&0&0\\
      \frac{1}{z-2}&1&0&0
    \end{bmatrix}&
    \mathbf{G}_2(z)&=
    \begin{bmatrix}
      1&0&0&0\\
      0&\frac{1}{z-2}&0&0\\
      0&0&0&0\\
      0&0&0&0
    \end{bmatrix}
  \end{align*}
  are network realizable on $\mathcal{G}$, but, as was pointed out in \cite{3994145}, this is not the case with their product
  \begin{align*}
    \mathbf{G}_1(z)\mathbf{G}_2(z)
    =\begin{bmatrix}
      0&0&0&0\\
      0&0&0&0\\
      \frac{1}{z-2}&\frac{1}{z-2}&0&0\\
      \frac{1}{z-2}&\frac{1}{z-2}&0&0
    \end{bmatrix}.
  \end{align*}
  This shows that the stability assumption is essential for statement \ii in Theorem~\ref{thm:add}.
\end{rmk}

\begin{pf*}{Proof of Theorem~\ref{thm:add}}
Let $\mathbf{G}_i(z)=C^i(zI-A^i)^{-1}B^i+D^i$ where
\begin{align*}
    \left[\begin{array}{c|c}
      A^i&B^i\\\hline C^i&D^i
    \end{array}\right]
    &=\left[\begin{array}{ccc|ccc}
      A^i_{11}&\ldots&A^i_{1N}&B^i_{1}&&0\\
      \vdots&&\vdots&&\ddots\\
      A^i_{N1}&\ldots&A^i_{NN}&0&&B^i_{N}\\\hline
      C^i_{11}&\ldots&C^i_{1n}&D^i_{1}&&0\\
      \vdots&&\vdots&&\ddots\\
      C^i_{N1}&\ldots&C^i_{NN}&0&&D^i_{N}
    \end{array}\right],
\end{align*}
Then $\mathbf{G}_3(z)=\mathbf{G}_1(z)+\mathbf{G}_2(z)$ provided that 
\begin{align*}
  A^3_{kl}&=\begin{bmatrix}A^1_{kl}&0\\0&A^2_{kl}\end{bmatrix}&
  B^3_{kl}&=\begin{bmatrix}B^1_{k}\\B^2_{k}\end{bmatrix}\\
  C^3_{kl}&=\begin{bmatrix}C^1_{kl}&C^2_{kl}\end{bmatrix}
\end{align*}
for all $k$ and $l$. The sparsity conditions, as well as stabilizability and detectability, follow trivially and \itemi holds.

Similarly, $\mathbf{G}_3(z)=\mathbf{G}_2(z)\mathbf{G}_1(z)$ provided that 
\begin{align}
  \left[\begin{array}{c|c}
    A^3_{kl}&B^3_{k}\\\hline
    C^3_{kl}&D^3_{k}
  \end{array}\right]
  &=\left[\begin{array}{cc|c}
    A^1_{kl}&0&B^1_{k}\\
    B^2_{k}C^1_{kl}&A^2_{kl}&B^2_{k}D^1_{k}\\\hline
    D^2_{k}C^1_{kl}&C^2_{kl}&D^2_{k}D^1_{k}
  \end{array}\right],
\label{eqn:product}
\end{align}
so $\mathbf{G}_3$ satisfies the sparsity conditions. Both factors are assumed to be stable, so stabilizability and detectability hold trivially. Hence \ii follows.

If $\mathbf{G}_1(z)=C(zI-A)^{-1}B+D$ and $D$ is invertible, then $\mathbf{G}_1^{-1}$ has the realization 
\begin{align*}
  \left[\begin{array}{c|c}
    A-BD^{-1}C&BD^{-1}\\\hline -D^{-1}C&D^{-1}
  \end{array}\right].
\end{align*}
where the needed sparsity structure, as well as stabilizability and detectability, follow from network realizability of $\mathbf{G}_1$. This proves \iii.
\end{pf*}

\section{Network Realizable Controllers}

The following theorem shows that in a number of cases, network realizability conditions on the controller can be mapped into similar conditions on closed loop transfer functions. From Theorem~1, we know that such constraints are convex, so they can be conveniently included in synthesis procedures based on convex optimization.

\begin{thm}
  Consider $\mathbf{P}$ and $\mathbf{C}$ such that $\mathbf{P}$ is strictly proper and define the closed loop matrix 
  \begin{align*}
    \mathbf{H}=\begin{bmatrix}
      I&-\mathbf{P}\\
      \mathbf{C}&I
    \end{bmatrix}^{-1}=\begin{bmatrix}
          (I+\mathbf{P}\mathbf{C})^{-1}&
          \mathbf{P}(I+\mathbf{C}\mathbf{P})^{-1}\\
          -\mathbf{C}(I+\mathbf{P}\mathbf{C})^{-1}&
          (I+\mathbf{C}\mathbf{P})^{-1}
        \end{bmatrix}.
  \end{align*}
  Then the following two statements are equivalent:
  \begin{description}
    \item[\itemi]Both $\mathbf{P}$ and $\mathbf{C}$ are network realizable on $\mathcal{G}$.
    \item[\ii]$\mathbf{H}$ is network realizable on $\mathcal{G}$.
  \end{description}
  The following two statements are also equivalent:
  \begin{description}
    \item[\iii]Both $\mathbf{P}\mathbf{C}$ and $\mathbf{C}$ are network realizable on $\mathcal{G}$.
    \item[\iv]Both $(I+\mathbf{P}\mathbf{C})^{-1}$ and $\mathbf{C}(I+\mathbf{P}\mathbf{C})^{-1}$ are network realizable on $\mathcal{G}$.
  \end{description}
  Suppose in addition that $\mathbf{P}$ is stable and network realizable on $\mathcal{G}$, while $\mathbf{H}$ is stable. Then the following are equivalent:
  \begin{description}
    \item[\itemv]$\mathbf{C}$ is network realizable on $\mathcal{G}$.
    \item[\vi]$\mathbf{C}(I+\mathbf{P}\mathbf{C})^{-1}$ is network realizable on $\mathcal{G}$.
  \end{description}
\label{thm:iff}
\end{thm}

\begin{rmk}
  The presentation in \cite{wang2016system} is focusing on finite impulse response representations of the closed loop maps. However, nothing excludes the use of other denominators when finite-dimensional parametrizations of closed loop dynamics are needed for computations. In fact, there is a rich literature on heuristics for selection of closed loop poles. 
\end{rmk}

\begin{rmk}
  The statement and proof of Theorem~\ref{thm:iff} is completely independent of how the set of stabilizing controllers is parametrized. The Youla-Kucera parametrization is the most well known option, but the parametrization suggested in \cite{wang2016system} appears to give simpler formulas for unstable plants.
\end{rmk}

\begin{rmk}
  Stability of the closed loop transfer matrix $\mathbf{H}$ means that all poles should be strictly inside the left half plane. In most applications the poles can actually restricted to a smaller subset $\Omega$ of the complex plane. Such a stronger assumption can be used to also get a stronger conclusion, namely that the controller has a network compatible realization with no uncontrollable or unobservable modes corresponding to poles outside $\Omega$.
\end{rmk}

\begin{pf*}{Proof of Theorem~\ref{thm:iff}.}
The strict properness of $\mathbf{P}$ implies that $\mathbf{H}(0)$ is invertible, so the equivalence between \itemi and \ii follows immediately from Theorem~\ref{thm:add}.

The equivalence between \iii and \iv follows from statement \iii in Theorem~\ref{thm:add} and the identity
\begin{align*}
  \begin{bmatrix}
    I+\mathbf{P}\mathbf{C}&0\\
    \mathbf{C}&I
  \end{bmatrix}^{-1}=\begin{bmatrix}
      (I+\mathbf{P}\mathbf{C})^{-1}&0\\
      -\mathbf{C}(I+\mathbf{P}\mathbf{C})^{-1}&I
    \end{bmatrix},
\end{align*}
since the strict properness of $\mathbf{P}$ gives both matrices an invertible direct term.

Assume that \itemv holds. Then \iii follows and therefore also \iv\!\!.  
This proves \vi\!\!.

Conversely, suppose that \vi holds. Then  $\mathbf{C}(I+\mathbf{P}\mathbf{C})^{-1}$ and $\mathbf{P}$ are both stable and network realizable on $\mathcal{G}$, so the same holds for their product $\mathbf{P}\mathbf{C}(I+\mathbf{P}\mathbf{C})^{-1}$. The identity $(I+\mathbf{P}\mathbf{C})^{-1}=I-\mathbf{P}\mathbf{C}(I+\mathbf{P}\mathbf{C})^{-1}$ gives stability and network realizable on $\mathcal{G}$ for $(I+\mathbf{P}\mathbf{C})^{-1}$. Hence \iv holds and the equivalence with \iii proves that $\mathbf{C}$ is network realizable on $\mathcal{G}$, so the proof is complete.
\end{pf*}

A common situation in applications is that a network realizable $\mathbf{Q}=\mathbf{C}(I+\mathbf{P}\mathbf{C})^{-1}$ has been designed and a corresponding controller needed. The equivalence between \itemv and \vi proves existence, but the proof of Theorem~\ref{thm:iff} is not convenient for construction of a corresponding controller. Instead, as will be seen in the next section, the classical Internal Model Control \cite{garcia1982internal} approach is useful for this task.

\section{Internal Model Control on Networks}
\begin{figure}
  \centerline{\includegraphics[width=.8\hsize]{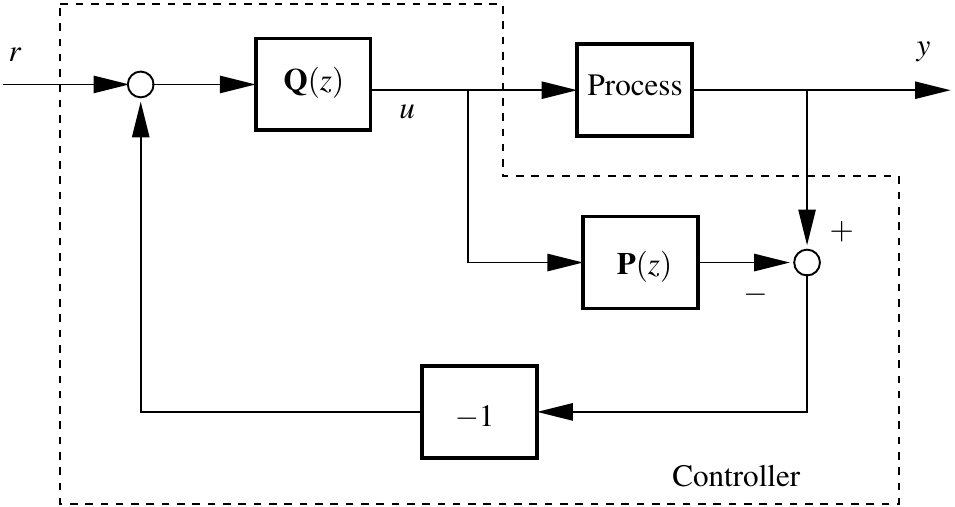}}
\caption{If $\mathbf{P}$ and $\mathbf{Q}$ are networks realizable, then a network realization of the controller can be obtained using Internal Model Control.}
\label{fig:IMC}
\end{figure}

Given $\mathbf{P}$ and $\mathbf{Q}$, consider a map from process output $y$ and reference value $r$ to control input $u$, defined by the equation
\begin{align}
  u=\mathbf{Q}[r+\mathbf{P}u-y].
\label{eqn:IMC}
\end{align}
See Figure~\ref{fig:IMC}.
Here $\mathbf{P}$ is the ``internal model'' that is used by the controller to predict the measured process output. The difference $\mathbf{P}u-y$ denotes a comparison between the predicted output $\mathbf{P}u$ and the measurement $y$. In the ideal case that the difference is zero, the control law reduces to $u=\mathbf{Q}r$, so $\mathbf{Q}$ defines the desired map from reference to input and $\mathbf{P}\mathbf{Q}$ is the resulting map from reference $r$ to output $y$. The transfer matrix from $r-y$ to $u$, given by (\ref{eqn:IMC}), is $\mathbf{C}=\mathbf{Q}(I-\mathbf{P}\mathbf{Q})^{-1}$.

Our interest in Internal Model Control stems from the fact it generates network realizable controllers in a very natural manner. Suppose that 
\begin{align*}
  \mathbf{P}(z)&=C(zI-A)^{-1}B\\
  \mathbf{Q}(z)&=G(zI-E)^{-1}F+H
\end{align*}
Then the controller $\mathbf{C}=\mathbf{Q}(I-\mathbf{P}\mathbf{Q})^{-1}$ has the realization
\begin{align*}
  \begin{bmatrix}\hat{x}^+\\\xi^+\end{bmatrix}
  &=\begin{bmatrix}A+BHC&BG\\FC&E\end{bmatrix}
  \begin{bmatrix}\hat{x}\\\xi\end{bmatrix}
  +\begin{bmatrix}BH\\F\end{bmatrix}(r-y)\\
  u&=\begin{bmatrix}HC&G\end{bmatrix}+H(r-y).
\end{align*}
It is easy to see that if $B$, $F$ and $H$ are diagonal, while $A$, $C$, $E$ and $G$ have a sparsity structure compatible with the graph $\mathcal{G}$, this realization is compatible with $\mathcal{G}$ after proper ordering of the states and block partitioning of the matrices. More specifically, if node $i$ of the graph is hosting the process state $x_i$, it should also host the controller state $(\hat{x}_i,\xi_i)$.

\begin{ex}
Consider a simple model for control of water levels in dams along a river:
\begin{align*}
\begin{cases}
  x_1(t+1)=0.9x_1(t)-u_1(t)\\
  x_2(t+1)=0.1x_1(t)+0.8x_2(t)+u_1(t)-u_2(t)\\
  x_3(t+1)=0.2x_2(t)+0.7x_3(t)+u_2(t)-u_3(t)
  \end{cases}
\end{align*}
Each state represents the water level in a dam and the control variables are used to control the release of water from one dam to the next. In this case, the transfer function from $(u_1,u_2,u_3)$ to $(x_1,x_2,x_3)$ is $\mathbf{P}(z)=(zI-A)^{-1}B$ where 
\begin{align*}
  A&=\begin{bmatrix}
    0.9&0&0\\
    0.1&0.8&0\\
    0&0.2&0.7
  \end{bmatrix}&
  B&=\begin{bmatrix}
    -1&0&0\\
    1&-1&0\\
    0&1&-1
  \end{bmatrix}
\end{align*}
A graph $\mathcal{G}$, corresponding to downwards flow of information, is defined by the node set $\mathcal{V}=\{1,2,3\}$ and the link set $\mathcal{E}=\{(1,1),(2,1),(2,2),(3,2),(3,3)\}$. The realization above does not have diagonal $B$-matrix, but the transfer function from $u$ to $x$ is still network realizable on $\mathcal{G}$, as shown by the (non-minimal) realization
\begin{align*}
  \bar{A}&={\small\left[\begin{array}{cc|cc|c}
    0.9&0&0&0&0\\
    0&0.8&0&0&0\\\hline
    0.1&0&0.8&0&0\\
    0&0.2&0&0.7&0\\\hline
    0&0&0.2&0&0.7
  \end{array}\right]}&
  \bar{B}&={\small\left[\begin{array}{r|r|r}
    -1&0&0\\
    1&0&0\\\hline
    0&-1&0\\
    0&1&0\\\hline
    0&0&-1
  \end{array}\right]}\\
  \bar{C}&={\small\left[\begin{array}{cc|cc|c}
    1&0&0&0&0\\\hline
    0&1&1&0&0\\\hline
    0&0&0&1&1
  \end{array}\right]}&
  \bar{D}&=0.
\end{align*}
Theorem~\ref{thm:iff} tells us that in order to find river dam controllers that only exchange information along the graph $\mathcal{G}$, it sufficient to consider $\mathbf{Q}=(I+\mathbf{C}\mathbf{P})^{-1}\mathbf{C}$ that are network realizable on $\mathcal{G}$. In particular, let $\mathbf{Q}$ have the state realization
\begin{align*}
  \left[\begin{array}{c|c}
    E&F\\\hline G&H
  \end{array}\right]
  &=\left[\begin{array}{ccc|ccc}
    E_{11}&0&0&F_{1}&0&0\\
    E_{21}&E_{22}&0&0&F_2&0\\
    0&E_{32}&E_{33}&0&0&F_3\\\hline
    G_{1}&0&0&0&0&0\\
    0&G_{2}&0&0&0&0\\
    0&0&G_{3}&0&0&0
  \end{array}\right]
\end{align*}
Then $\mathbf{C}=(I-\mathbf{Q}\mathbf{P})^{-1}\mathbf{Q}$ mapping $(e_1,e_2,e_3)$ to $(u_1,u_2,u_3)$ can be implemented as the Internal Model Controller
\begin{align*}{\footnotesize
  \!\!\!\begin{bmatrix}
  \hat{x}_1^+\\\xi_1^+\\\hat{x}_2^+\\\xi_2^+\\\hat{x}_3^+\\\xi_3^+  
  \end{bmatrix}}
  &=\left[{\footnotesize\begin{array}{cc|cc|cc}
      0.9&-G_{1}&0  &  0&  0& 0\\
      F_1  &E_{11}&0  &  0&  0& 0\\\hline
      0.1&G_{1}&0.8&-G_2&  0& 0\\
      0  &E_{21}&F_2  &E_{22}&  0& 0\\\hline
      0  &0  &0.2& G_2&0.7&-G_3\\
      0  &0  &0  &E_{32}&F_3  &E_{33}
    \end{array}}\right]
{\footnotesize\begin{bmatrix}
\hat{x}_1\\\xi_1\\\hat{x}_2\\\xi_2\\\hat{x}_3\\\xi_3
\end{bmatrix}
-\begin{bmatrix}
0\\e_1\\0\\e_2\\0\\e_3
\end{bmatrix}}
\end{align*}
This realization has a structure compatible with the graph $\mathcal{G}$ (in spite of the fact that it is based on the original state matrices $A,B$ rather than $\bar{A},\bar{B}$).
\end{ex}

\section{Acknowledgement}
Financial support from the Swedish Research Council and the Swedish Foundation for Strategic Research is gratefully acknowledged, The author is a member of the Linnaeus
center LCCC and the excellence center ELLIIT.

\end{document}